# A note on using Bayes priors for Weibull distribution

Pasquale Erto, *Member*, *IEEE*, and Massimiliano Giorgio

**Abstract—** In this note, the practical use of priors for Bayes estimators – of the two parameters of the Weibull reliability model – is discussed in a technological context. The meaning of the priors as expression of *virtual* data samples is analyzed. The implications of physics of failures are also highlighted. The whole analysis shows a *rational* way to convert a really available technological knowledge into prior information, effectively and without a long elicitation process. A large Monte Carlo study, on both complete and censored samples of small size, shows the good properties of some estimators, which can exploit such a new approach.

*Index Terms*—Bayes methods, Lifetime estimation, Parameter estimation, Reliability engineering, Weibull distribution,

## NOTATIONS

| | |
|---|---|
| $x$ | life length |
| $Sf(x)$ | reliability or survival function, $Sf(x) = \Pr\{X > x\}$ |
| $pdf(x)$ | probability density function |
| $n$ | sample size |
| $\boldsymbol{x}$ | $n$ dimensional type II censored random sample |
| $r$ | number of failure times in the random sample $\boldsymbol{x}$ |
| $x_i$ | indicates the $i$-th life length |
| $D$ | indicates the set of the subscripts of the observed death (failure) times |

Pasquale Erto is with the Industrial Engineering Department, University of Naples Federico II, Naples 80125, Italy (e-mail: ertopa@unina.it).

Massimiliano Giorgio is with the Industrial and Information Engineering Department, Second University of Naples, Aversa 81031, Italy (e-mail: massimiliano.giorgio@unina2.it ).

| | |
|---|---|
| $C$ | indicates the set of subscripts of censored lifetimes |
| $\alpha, \beta$ | scale and shape parameters of the Weibull distribution |
| $R$ | reliability level for the referred reliable life |
| $x_R$ | reliable life, viz. quantile of the distribution such as that $Sf(x_R) = R$ $(0 \leq R \leq 1)$ |
| $\bar{x}_R$ | prior numerical value anticipated for $x_R$ |
| $\beta_1, \beta_2$ | limits of the prior numerical interval anticipated for $\beta$ |
| $\wedge$ | implies a ML estimator |
| $\sim$ | implies a PBE estimator |
| $B_{n,r}$ | indicates the unbiasing factor such that $E[B_{n,r} \hat{\beta}] = \beta$, viz. $B_{n,r} = 1/E[\hat{\beta}/\beta]$). This factor depends only on $n$ and $r$ (see as instance [1]) |
| $\bar{\beta}$ | indicates the unbiased MLE of $\beta$, $\bar{\beta} = B_{n,r} \hat{\beta}$ |
| $DS[\cdot]$ | indicates the standard deviation of the estimator in parenthesis |
| $RQ[\cdot]$ | indicates the square root of the mean square error of the estimator in parenthesis, viz. $RQ[\cdot] = \left(DS^2[\cdot] + \Delta^2[\cdot]\right)^{0.5}$ |
| $\Gamma(z)$ | is the Gamma function, viz. $\Gamma(z) = \int_0^\infty u^{z-1} e^{-u} du, \quad \forall z > 0$ |

I. INTRODUCTION

Typically, who uses Bayesian methods to estimate distribution parameters, chooses this approach because he possesses other quantifiable information besides a sample of experimental data. Unfortunately, practical use of such information is often associated with difficulties related to its elicitation and its formalization into prior distributions. However, in technological field, these difficulties can be faced more easily. In fact, the technological systems, by their nature, are the products of engineers who attain precise knowledge about *some their characteristics that are more or less explicitly related to* the parameters of their reliability model.

On the other hand, the high rhythm of product innovation and the high cost per unit of the technological systems/components often strongly limit the number of experimental data (say to 3-5) available for the analysis. It is not difficult to perceive that in these particular experimental situations the effectiveness of classical statistical methods is drastically reduced. In such conditions, classical methods may supply estimates that are even worse than those which can be anticipated on the basis of elementary technological considerations. Moreover, even if "reasonable" estimates were obtained, it would be however difficult to give evidence of their "objectivity". It doesn't exist a test that, on the basis of so few data, can demonstrate the goodness of the method and/or of the model used for the analyses.

Therefore, in the above cases, being compelled to operate a subjective choice, we have a further reason to use a Bayesian approach. In fact, in its framework we are obliged to do "subjective choices" in more conscious and explicit way. Often, the question we have to reply to is not "what can we conclude on the basis of these data?" but rather "what can we decide using these data too?"

Apart of the Reliability, also the Statistical Process Control is involved in the above considerations. In fact, when a key quantity associated with the life length of whichever thing must be monitored, the solution is a control chart of the key quantity, whose variability is often modeled via skewed distribution. Unfortunately, classical Shewhart control charts stand under *s*-normality assumption. So they are not effective when the distribution of the key quantity is skewed and the small size of the available samples prevents the exploitation of the "*s*-normalizing effect". In such a case, the recourse to Bayesian methods appears appropriate. So, in the recent literature Bayesian methods that combine *initial* prior information with *current* sample data are being taking into growing consideration [5] [6] [11] [12] [23] - [26].

## II. BAYESIAN ESTIMATORS FOR THE WEIBULL RELIABILITY MODEL

The Weibull distribution is one of the most commonly used distributions in technological field. It does not have therefore to amaze that the problem of the estimation of the parameters of this model has been treated with particular attention by Bayesian literature. The most interesting case for the applications is surely the one in which both its scale and the shape parameters have to be estimated. Unfortunately, as asserted by Soland [30], for the two-parameter Weibull model a conjugate family of continuous joint prior distributions doesn't exist. This situation has played in the Bayesian context the same role played in the classical one by the lack of joint sufficient statis-tics of fixed dimension (e.g., see [16]). That produced from the 70's to the 90's a strong stimulus to formulate Bayesian estimators, each proposing the use of even more than one type of prior distribution (for a critical analysis see [19]). As an example, the use of mixed prior distributions (discrete for the shape parameter, continuous for scale parameter) has been suggested in [30]; the use of continuous prior distributions (Uniform for the shape parameter and Inverted Gamma for the scale parameter) has been proposed in [31]; many different prior distributions both discrete and continuous have been proposed in [22] (Inverted Gamma - Compound Inverted Gamma, Discrete mass function - Compound Inverted Gamma, as well as Uniform distribution - Compound Inverted Gamma, respectively for the shape and scale parameter). In [15] a gamma prior on scale parameter and no specific prior on shape parameter is assumed (i.e., it is only assumed that the support of the shape parameter is $0, \infty$ and its density function is of log-concave type). The gamma prior on both the scale and shape parameters are considered in [3]. Although of a great interest, also the last two approaches may be challenging to apply to real life problems due to the difficulties of evaluating the needed prior information.

Then, the majority of the methods proposed in these papers can be object of criticism, having the characteristic of considering prior distributions selected more for their tractability than for the effective ability to represent the available prior information. Moreover, these distributions often refer to parameters whose meaning is completely unknown to technologists. That forces the user to operate a transformation/transfer of the prior information he possess, from the domain in which it effec-tively stands to that one in which it must be formulated and then used.

From the above discussion, the following general conclusions can be schematically drawn:

I. It is a fact that in Reliability it exists a large demand for Bayesian methods of estimation of the Weibull model parameters that are easy to implement (as an example, methods that do not need the use laborious techniques or that allow a natural elicitation of the prior information). Also the Statistical Process Control is involved in this demand, since often the key quantity to be monitored shows a variability modeled via Weibull distribution (e.g. see [20]).

II. Technicians/design engineers are nearly always in a position to say whether the mechanism of failure that attempts to their product depends on time or not, and if it depends, whether that involves an improvement or a worsening. Therefore they can establish, theoretically, if the shape parameter of the Weibull is smaller than 1 (improvement), greater than 1 (worsening) or equal to 1 (independence). Moreover the technicians often know the results of many previous experiences that allow them to statistically define, an interval, of all equally plausible values, that contains the shape parameter with a high level of confidence. The typically defined intervals are $[0.5, 1]$, $[1, 3]$, $[0.5, 2]$, but even in the case of a wider interval, like $(0, 10]$ (that contains the totality of the estimates of the shape parameter obtained in the history of technology) the technician would however enable the Bayesian method to spend the experimental information on an interval infinitely smaller than the shape parameter range $(0, \infty)$.

III. Nearly all the projects of technological systems contain (not statistical) assessments of reliable life or quantile $x_R$ (where $R$ is a prefixed reliability level) being this a required technical specification. It is therefore natural to think that in such cases it is possible to anticipate an a priori estimate of this parameter. Obviously, a remarkable level of uncertainty is associated to this kind of estimate. This uncertainty is *inversely proportional to the expertise* of the technician that supplies it. Besides, also in Statistical Process Control a technical specification in terms of Weibull percentile (e.g., as a minimum threshold for reliability design) must be very often monitored instead of other parameters such as the scale parameter or the mean [2]. So, in the recent literature the need to monitor a process percentile under Weibull assumption is arisen [5] [6] [11] [12].

IV. The meaning of the prior distribution parameters is extraneous to the technicians' mentality as well as to the their cultural heritage and/or experiences.

On the basis of these considerations some attempts were made in [8] - [10] to allow a genuine and direct use (i.e., in the same form in which they are commonly possessed) of the prior information that is typically available in technological application. More recently, a novel procedure has been suggested in [14] where the use of prior information in the form of the interval assessment of the reliability function (as opposed to that on the Weibull parameters) is proposed. When such specific information is available, this procedure allows constructing continuous joint prior distribution of Weibull parameters very effectively.

The present paper starts from the results in [8] - [10] and presents a new *rational*

approach to define prior distribution parameters that allow better fitting the prior uncertainty.

## III. GENERALIZED PRIORS

### A. Assumptions

1. The considered random variable, $X$, has a two-parameters Weibull reliability function:

$$Sf(x|\alpha,\beta) = \exp\left[-(x/\alpha)^\beta\right], \qquad x \geq 0; \quad \alpha,\beta > 0. \qquad (1)$$

2. The expert (that can eventually be the user himself) is in a position to anticipate the limits, $\beta_1$ and $\beta_2$ (where $\beta_2 > \beta_1 > 0$), of an interval of values in which the unknown shape parameter $\beta$ is contained and is uniformly distributed (point II of par. II).

3. The user, either on the basis of his own experience, asking the advise of an expert or utilizing empirical tables [27], is in a position to anticipate a value for the parameter $x_R$ (point III of par. II).

4. The user can specify the level of uncertainty of the value anticipated for the parameter $x_R$. Such uncertainty can be used to quantify the weight to give to the expert opinion on which the empirical estimate $\bar{x}_R$ is based (point III of par. II).

### B. Re-parameterization of the Weibull model

In order to facilitate the elicitation of the prior distributions, we will use the follow-ing re-formulation of the Weibull model:

$$Sf(x|x_R,\beta) = \exp\left[-K(x/x_R)^\beta\right]; \qquad x \geq 0; \quad x_R, \beta > 0; \quad K = \log(1/R) \qquad (2)$$

which is expressed in terms of the parameters, which the prior information is referred to (points 2, 3 and 4 of par. A).

## C. *The prior distribution of the shape parameter* $\beta$

As prior distribution for the parameter $\beta$ the following continuous Uniform distribution was adopted:

$$pdf(\beta) = 1/(\beta_2 - \beta_1); \qquad \beta_2 \geq \beta \geq \beta_1 > 0; \quad \beta_2 > \beta_1 \qquad (3)$$

that is able to describe in simple and not restrictive form the kind of information considered at point 2 of par. A.

It can be worth to observe that:

- The prior (3) is an informative prior distribution. In fact it hasn't been used in order to represent lack of information about the form of $pdf(\beta)$ in the interval $[\beta_1, \beta_2]$, but in order to effectively describe the kind of information, which is usually available. That implies a probability uniformly distributed on the interval $[\beta_1, \beta_2]$, which certainly includes $\beta$ on the basis of physics considerations and experience.

- Suitability of (3) cannot be evaluated on the basis of pure mathematical consideration. In fact, this approach would not take into account that the information regarding the shape parameter, $\beta$, are typically based on the study of physics of failure

(point II of par. II).

## D. *The prior distribution of the parameter* $x_R$

Considered the nature of the prior information (point 2, 3 and 4 of par. A) we adopted the following Inverted Generalized Gamma *pdf* as conditional (i.e., given $\beta$) prior of $x_R$ (see Appendix of [10] and Fig. 1):

$$pdf(x_R|\beta) = \frac{\beta \, a^{\beta w}}{\Gamma(w)} x_R^{-(\beta w + 1)} e^{-(x_R/a)^{-\beta}}; \quad a, w > 0. \tag{4}$$

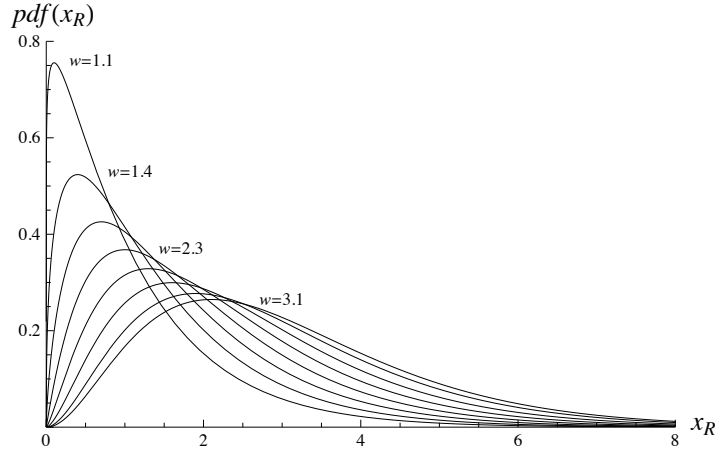

Fig. 1. Probability density function of the Inverted Generalized Gamma model, for $a = 1$, $\beta = 1$ and $w = 1.1\,(0.3)\,3.1$ (these settings are close to those used in par. IV).

Even though it has been emphasized that the existence of a prior can be always included in a system of axioms [28], in order to facilitate the practical use of this prior, we observe that:

- It goes to zero as $x_R$ goes to $0^+$, and has unique mode, which exists for any value of its parameters. Differently from many other prior models commonly used for the Weibull (for a deep analysis see [19]), it is not of negative exponential type as $x_R$ goes to infinity, but it is infinitesimal of lower order (i.e., $w\beta+1$). So, it is less restrictive since it reserves a heavy probability tail even to values much greater than those thought more probable.

- It is a natural conjugate (given $\beta$) of the Weibull distribution. This property allows assigning the appropriate weights to both the experimental observations and the prior information in a rational way. In fact, combining via Bayes theorem the prior (4) and the likelihood function of the random sample, $x$:

$$L(x_R,\beta|x) \propto \prod_{i \in C} Sf(x_i|x_R,\beta) \prod_{i \in D} pdf(x_i|x_R,\beta) =$$

$$= \left(\frac{K\beta}{x_R^\beta}\right)^r P^{\beta-1} \cdot \exp\left[-\frac{K\,S(\beta)}{x_R^\beta}\right]; \quad (5)$$

$$\text{where} \quad K = \ln(1/R), \quad S(\beta) = \sum_{i=1}^{n} x_i^\beta, \quad P = \prod_{i \in D} x_i$$

we obtain the following conditional posterior distribution:

$$pdf(x_R|\beta,x) = \frac{\beta\left[a^\beta + K\,S(\beta)\right]^{w+r}}{\Gamma(w+r)} x_R^{-(w+r)\beta-1} e^{-\left[a^\beta + K\,S(\beta)\right]x_R^{-\beta}}, \quad (6)$$

that coincides with (4) where $w$ and $a^\beta$ are replaced with $w+r$ and $a^\beta + K\,S(\beta)$ respectively (being $r$ the number of failure times in the random sample $x$ of size $n$) This precious computational feature allows also understanding that the role-played in (6) by $a^\beta$ and $w$ is equivalent to that played by $K\,S(\beta)$ and $r$ respectively.

- Adopting an approach similar to that proposed in [7] and [17], we can consider

(4) as at the posterior (conditional) distribution obtained combining the likelihood function of a virtual random sample $x'$ (representing the virtual source of the prior information) and the Jeffreys non-informative prior [4], $pdf(x_R|\beta) \propto x_R^{-1}$. This prior correctly expresses the state of "relative" ignorance in which the elicitators were before gaining the information supplied by the virtual sample $x'$. This approach leads to rewrite (4) as:

$$pdf(x_R|\beta, x') = \frac{\beta \left(K\, S'(\beta)\right)^{r'}}{\Gamma(r')} x_R^{-r'\beta-1} e^{-(K\, S'(\beta))x_R^{-\beta}}; \qquad S'(\beta) = \sum_{i=1}^{n'}(x'_i)^{\beta} \qquad (7)$$

where $k\, S'(\beta)$ and $r'$ stand as $a^\beta$ and $w$ in (4) (the prime " ′ " indicates quantity based on the virtual sample $x'$). This consideration and those about (6) can be used as powerful and practical tools to rationalize both the elicitation of the prior information and the calibration of prior distribution.

- The interpretative result (7) evidences that the prior *pdf* (4) allows managing prior information about $x_R$ and the experimental data in an equivalent manner. Indead, it gives us chance to look at (6) as a result obtained by combining the Jeffreys non-informative prior and the likelihood function of both real and virtual samples, $x$ and $x'$ (see also [7] and [17]). In practice, we effectively treat both the pieces of information homogeneously. In fact, first, we collect prior information about $x_R$ by using empirical and substantial statistical procedures (e.g.: averaging on past experiments; using information relative to similar devices; using empirical tables). Then, we use the real sampling data to refine the form of the prior distribution.

- The prior *pdf* model (4) includes many others already proposed in literature. For instance, setting:

$$w = 1 \qquad (8)$$

we obtain the distribution used in [8]. Instead, placing:

$$w = \begin{cases} 1, & \text{if } \beta \geq 1 \\ 1/\beta^2, & \text{if } \beta < 1 \end{cases} \qquad (9)$$

we obtain the prior used in [9].

## E. Converting prior information available for $x_R$

### 1) Setting the hyperparameter a

In order to set the hyperparameter $a$, the following relation can be used:

$$a = \overline{x}_R \frac{\Gamma(w)}{\Gamma(w - 1/\beta)}; \qquad w > 1/\beta. \qquad (10)$$

It is easy to verify that, being $w > 1/\beta$ and $\beta \in [\beta_1, \beta_2]$, by using (10) it results:

$$E[x_R | \beta] = \int_0^\infty x_R \, pdf(x_R | \beta) \, dx_R = \overline{x}_R; \qquad \beta > 0. \qquad (11)$$

Therefore, obviously, it results also:

$$E[x_R] = \int_{\beta_1}^{\beta_2} \int_0^\infty x_R \, pdf(x_R | \beta) \, pdf(\beta) \, dx_R \, d\beta = \overline{x}_R, \qquad (12)$$

regardless of the values given to the $\beta_1$ and $\beta_2$ parameters of the prior (3).

We can observe that (10) does not predetermine the $a$ value, but it only transforms the $a$ value in a function of $\beta$ and $w$. Moreover comparing (10) to the relation:

$$a^\beta = K\, S'(\beta) = K \sum_{i=1}^{n'} (x'_i)^\beta; \qquad K = \log(1/R), \qquad (13)$$

suggested by (7), it can be deduced that (10) allows obtaining, for every given value of $\beta$ and $w$, the only value of the statistics $K\, S'(\beta)$ that satisfies the relation (11).

In order to further clarify the motivations and the role of the adopted relationship (10), we observe:

- The indirect evaluation of $a$ is surely easier to implement and more effective than the direct evaluation by means of (13). In fact to directly use the relationship (13) we would be obliged to known the whole virtual random sample, rather than the single value $\bar{x}_R$.
- Using $\bar{x}_R$ as prior mean of the parameter $x_R$ is equivalent to consider $\bar{x}_R$ a prior Bayesian point estimate of $x_R$;
- Using (10) implies that the value of the parameter $\beta$ cannot affect the meaning of $\bar{x}_R$, since it can only modify the dispersion of the prior distribution of $x_R$ (i.e., the level of uncertainty associated to the prior empirical estimate $\bar{x}_R$).

*2) Setting the hyperparameter w*

By means of the value given to $w$, we can comparatively quantify the weights to give to prior information (point 4 of par. A) and experimental data. In fact, fixed the $a$ and $\beta$ values, the higher $w$ is the lower the dispersion of (4) is. This proportionality still stands even if $a$ is calculated by using (10) given $\bar{x}_R$.

Specifically, some suggestions to set the $w$ value can be obtained analyzing the *pdf* (6) and (7). In (6) the hyperparameters $w$ and $a$ and the experimental data are combined to obtain the parameters of the posterior *pdf* (6) conditioned to $\beta$:

$$w + r = r' + r \tag{14}$$

$$a^\beta + K\, S(\beta) = K\left[\sum_{i=1}^{n'}(x'_i)^\beta + \sum_{i=1}^{n} x_i^\beta \right], \tag{15}$$

(where the prime " ′ " indicates virtual quantity). Since (10), the higher $w$ is the higher the $w$ and $a^\beta$ contributions to the posterior parameters (14) and (15) are. So, modulating the value given to $w$, we can predetermine the weight given to prior information via posterior *pdf* parameters.

It is worth to remarks that, even if the consideration about (6) and (7) have been done in terms of integer value of $w$, from a mathematical point of view we can assign $w$ any positive value. However, if we agree to use (10), we must give it a value greater than $1/\beta$.

*F. The joint posterior distribution*

Combining, via Bayes theorem, the joint prior distribution:

$$pdf(x_R, \beta) = pdf(x_R|\beta)\, pdf(\beta) = \frac{1}{\beta_2 - \beta_1} \frac{\beta\, a^{\beta w}}{\Gamma(w)} x_R^{-w\beta - 1} \exp^{-(x_R/a)^{-\beta}}$$

and the likelihood function (5), the following joint posterior distribution is obtained:

$$pdf(x_R, \beta | \mathbf{x}) = \frac{\beta^{r+1} \, a^{\beta w} \, x_R^{-(r+w)\beta-1} \, P^\beta \, \exp[-x_R^{-\beta} A] \, \Gamma^{-1}(w)}{\int_{\beta_1}^{\beta_2} \beta^r \, a^{\beta \cdot w} \, P^\beta \, A^{-(w+r)} \, \Gamma(w+r) \, \Gamma^{-1}(w) \, d\beta} \tag{16}$$

where:

$$A = K \left( \sum_{i=1}^{n} x_i^\beta \right) + a^\beta, \qquad P = \prod_{i \in D} x_i, \qquad k = \log(1/R).$$

Using such a distribution, we can formulate point or interval estimators of any combination of the Weibull parameters (shape and scale, shape and reliable life etc.).

*1) Bayes estimators*

The point estimators, $\tilde{x}_R$ and $\tilde{\beta}$, of $x_R$ and $\beta$ can be obtained using the following relations:

$$\tilde{x}_R = E[x_R | \mathbf{x}] = I_1/I_0; \qquad \tilde{\beta} = E[\beta | \mathbf{x}] = I_2/I_0 \tag{17}$$

where:

$$I_h = \int_{\beta_1}^{\beta_2} \beta^{r_h} \, a^{\beta \cdot w} \, P^\beta \, A^{-[r+w-m_h]} \, \Gamma(r+w-m_h) \, \Gamma^{-1}(w) \, d\beta \tag{18}$$

being:

$$h = 0, 1, 2; \qquad r_0 = r_1 = r; \qquad r_2 = r+1; \qquad m_0 = m_2 = 0; \qquad m_1 = 1/\beta.$$

Since the integrals (18) cannot be solved analytically, Bayes estimators have to be performed numerically.

# IV. MONTE CARLO STUDY OF THE PERFORMANCES OF BAYES ESTIMATORS

The non asymptotic properties of the point Bayes estimators have been evaluated performing a Monte Carlo study based on the use of three groups of 2000 complete samples of size $n=3$ and three groups of 2000 type II censored samples with $n=5$ and $r=3$. These samples have been generated from three Weibull distributions with shape parameter $\beta$ equal to 2, 1 and 0.6 respectively, and the reliable life $x_R$ equal to 1 (with $R=0.98$) in all cases.

For every $\beta$ we considered nine possible different kind of prior information to represent as many as possible experimental situations. For the shape parameter $\beta$ we used prior intervals wider than those typically suggested in literature (e.g., see [21]).

The used prior information is shown in Tables 1 and 2. For every combination of prior information, we examined the performances of the Bayes estimators in terms of bias, standard deviation and root mean square error of their empirical distribution. The performances of the Bayes estimators have been compared with those of corresponding MLE.

We assigned many different values to the hyperparameter $w$, excluding those values that would have produced excessively strong (i.e., dominant, [4]) prior distributions. In this way we tried to leave the estimators "learning" from the considered very small data sample. Specifically, we set $w=1.1/\beta$, $w=1.4/\beta$, $w=1.8/\beta$ and $w=1/\beta_1+0.1$ in order to both satisfy the constraint of (10) and guarantee values smaller than the number of failure, $r$ (see (14)). Moreover, setting $w=c/\beta$ guarantees large dispersion of the prior (4) for any value of $\beta$, because the *pdf* (4) results to be infinitesimal of order $c+1 = w\beta +1$ (see par. D) as $x_R$ goes to infinity and $\beta \in [\beta_1, \beta_2]$. Obviously, as noted before, the above $w$ setting does not fix the hyperparameter $w$ but transform it in a function of $\beta$.

## A. Prior information used in simulation

Type 1 intervals used for $\beta$ in Table 1 are centered on the true value; the Type 2 are upper biased, that is the true value is equal to the lower bound of the intervals; the Type 3 intervals are lower biased, that is the true value of $\beta$ is equal to the upper bound of the intervals (the true value is that used to generate the pseudo random samples).

For $\bar{x}_R$ in Table 2 we used values 1, 10 and 0.1, respectively equal to the true value, ten time greater than the true value and ten time smaller than the true value (the true value is that used to generate the pseudo random samples)

Table 1. Prior interval used in simulation for the parameter $\beta$.

| Type # | Prior interval of $\beta$ | | |
|---|---|---|---|
| | ($\beta$=2) | ($\beta$=1) | ($\beta$=0.6) |
| 1 | 1-3 | 0.7-1.3 | 0.3-0.9 |
| 2 | 2-4 | 1-1.3 | 0.6-0.9 |
| 3 | 0.5-2 | 0.7-1 | 0.3-0.6 |

Table 2. Combinations of prior information for $x_R$ and $\beta$ used in the simulation. The true value of $x_R$ is always 1.

|  |  | \multicolumn{3}{c}{Prior values of $\bar{x}_R$} |  |  |
|---|---|---|---|---|
|  | Type # | 1 | 10 | 0.1 |
| **Prior interval used for $\beta$** | 1 | I | II | III |
|  | 2 | IV | V | VI |
|  | 3 | VII | VIII | IX |

## V. RESULTS

All the obtained results of the simulation study are presented in the following Tables 3-8. To understand these Tables both the Table 1 and Table 2 information must be taken into account. For instance the Test "VI" refers to the notation used in Table 2 and means the use of:

- a lower biased value for $x_R$ (i.e., a value ten times smaller then the "true" value $x_R = 1$);

- a Type 2 biased interval for $\beta$ (i.e., the interval 2, 4 if the "true" value of $\beta$ is 2; the interval 1, 1.3 if the "true" value of $\beta$ is 1; the interval 0.6, 0.9 if the "true" value of $\beta$ is 0.6).

Table 3. Performances of the Bayes estimators for $\beta = 2$, $n = r = 3$ and different $w$ values.

| $RQ[\tilde{x}_R]$ | | | | | $RQ[\tilde{\beta}]$ | | | |
|---|---|---|---|---|---|---|---|---|
| $w$ values | | | | | $w$ values | | | |
| $1.1/\beta$ | $1.4/\beta$ | $1.8/\beta$ | $1/\beta_1+0.1$ | Test # | $1.1/\beta$ | $1.4/\beta$ | $1.8/\beta$ | $1/\beta_1+0.1$ |
| .38E+00 | .29E+00 | .23E+00 | .23E+00 | I | .34E+00 | .25E+00 | .22E+00 | .22E+00 |
| .41E+00 | .12E+01 | .21E+01 | .13E+01 | II | .28E+00 | .48E+00 | .44E+00 | .55E+00 |
| .44E+00 | .51E+00 | .60E+00 | .61E+00 | III | .42E+00 | .50E+00 | .60E+00 | .59E+00 |
| .88E+00 | .82E+00 | .75E+00 | .85E+00 | IV | .75E+00 | .71E+00 | .68E+00 | .79E+00 |
| .97E+00 | .16E+01 | .25E+01 | .14E+01 | V | .86E+00 | .12E+01 | .12E+01 | .52E+00 |
| .88E+00 | .81E+00 | .72E+00 | .61E+00 | VI | .75E+00 | .70E+00 | .63E+00 | .50E+00 |
| .63E+00 | .48E+00 | .40E+00 | .28E+00 | VII | .71E+00 | .53E+00 | .45E+00 | .38E+00 |
| .21E+00 | .74E+00 | .16E+01 | .20E+01 | VIII | .35E+00 | .33E+00 | .40E+00 | .91E+00 |
| .82E+00 | .81E+00 | .82E+00 | .87E+00 | IX | .10E+01 | .96E+00 | .96E+00 | .10E+01 |

Table 3b. Performances of the MLE for $\beta = 2$ and complete sampling.

| $n=r$ | $RQ[\hat{x}_R]$ | $RQ[\hat{\beta}]$ | $DS[\bar{\beta}]$ |
|---|---|---|---|
| 3 | .21E+01 | .75E+01 | .31E+01 |
| 5 | .13E+01 | .18E+01 | .11E+01 |
| 7 | .98E+00 | .11E+01 | .78E+00 |
| 10 | .74E+00 | .77E+00 | .59E+00 |
| 15 | .56E+00 | .54E+00 | .45E+00 |
| 22 | .43E+00 | .40E+00 | .36E+00 |
| 30 | .35E+00 | .33E+00 | .30E+00 |

Table 4. Performances of the Bayes estimators for $\beta = 1$, $n = r = 3$ and different $w$ values.

| $RQ[\tilde{x}_R]$ | | | | | $RQ[\tilde{\beta}]$ | | | |
|---|---|---|---|---|---|---|---|---|
| $w$ values | | | | | $w$ values | | | |
| $1.1/\beta$ | $1.4/\beta$ | $1.8/\beta$ | $1/\beta_1+0.1$ | Test # | $1.1/\beta$ | $1.4/\beta$ | $1.8/\beta$ | $1/\beta_1+0.1$ |
| .47E+00 | .38E+00 | .32E+00 | .41E+00 | I | .97E-01 | .76E-01 | .70E-01 | .57E-01 |
| .71E+00 | .16E+01 | .24E+01 | .17E+01 | II | .60E-01 | .10E+00 | .12E+00 | .61E-01 |
| .53E+00 | .57E+00 | .64E+00 | .53E+00 | III | .13E+00 | .15E+00 | .17E+00 | .12E+00 |
| .11E+01 | .10E+01 | .89E+00 | .12E+01 | IV | .12E+00 | .13E+00 | .13E+00 | .15E+00 |
| .13E+01 | .19E+01 | .26E+01 | .17E+01 | V | .15E+00 | .17E+00 | .17E+00 | .15E+00 |
| .11E+01 | .94E+00 | .77E+00 | .10E+01 | VI | .12E+00 | .11E+00 | .10E+00 | .13E+00 |
| .58E+00 | .49E+00 | .42E+00 | .50E+00 | VII | .17E+00 | .15E+00 | .14E+00 | .14E+00 |
| .32E+00 | .82E+00 | .17E+01 | .53E+00 | VIII | .13E+00 | .11E+00 | .11E+00 | .15E+00 |
| .65E+00 | .68E+00 | .72E+00 | .63E+00 | IX | .19E+00 | .19E+00 | .20E+00 | .17E+00 |

Table 4b. Performances of the MLE for $\beta = 1$ and complete sampling.

| $n=r$ | $RQ[\hat{x}_R]$ | $RQ[\hat{\beta}]$ | $DS[\bar{\beta}]$ |
|---|---|---|---|
| 3 | .13E+02 | .37E+01 | .16E+01 |
| 5 | .63E+01 | .91E+00 | .56E+00 |
| 7 | .41E+01 | .56E+00 | .39E+00 |
| 10 | .26E+01 | .38E+00 | .30E+00 |
| 15 | .17E+01 | .27E+00 | .23E+00 |
| 22 | .12E+01 | .20E+00 | .18E+00 |
| 30 | .90E+00 | .16E+00 | .15E+00 |

Table 5. Performances of the Bayes estimators for $\beta = 0.6$, $n = r = 3$ and different $w$ values.

| $RQ[\tilde{x}_R]$ | | | | | $RQ[\tilde{\beta}]$ | | | |
|---|---|---|---|---|---|---|---|---|
| $w$ values | | | | | $w$ values | | | |
| $1.1/\beta$ | $1.4/\beta$ | $1.8/\beta$ | $1/\beta_1+0.1$ | Test # | $1.1/\beta$ | $1.4/\beta$ | $1.8/\beta$ | $1/\beta_1+0.1$ |
| .53E+00 | .46E+00 | .37E+00 | .28E+00 | I | .11E+00 | .78E-01 | .74E-01 | .82E-01 |
| .19E+01 | .36E+01 | .46E+01 | .50E+01 | II | .72E-01 | .11E+00 | .13E+00 | .12E+00 |
| .82E+00 | .81E+00 | .84E+00 | .87E+00 | III | .18E+00 | .16E+00 | .15E+00 | .16E+00 |
| .36E+01 | .28E+01 | .21E+01 | .36E+01 | IV | .98E-01 | .10E+00 | .10E+00 | .12E+00 |
| .46E+01 | .53E+01 | .58E+01 | .57E+01 | V | .13E+00 | .15E+00 | .16E+00 | .15E+00 |
| .33E+01 | .23E+01 | .14E+01 | .28E+01 | VI | .87E-01 | .77E-01 | .65E-01 | .95E-01 |
| .77E+00 | .61E+00 | .51E+00 | .56E+00 | VII | .15E+00 | .12E+00 | .10E+00 | .12E+00 |
| .41E+00 | .91E+00 | .18E+01 | .15E+01 | VIII | .10E+00 | .83E-01 | .77E-01 | .11E+00 |
| .92E+00 | .90E+00 | .89E+00 | .90E+00 | IX | .20E+00 | .17E+00 | .16E+00 | .17E+00 |

Table 5b. Performances of the MLE for $\beta = 0.6$ and complete sampling.

| $n=r$ | $RQ[\hat{x}_R]$ | $RQ[\hat{\beta}]$ | $DS[\bar{\beta}]$ |
|---|---|---|---|
| 3 | .17E+03 | .22E+01 | .94E+00 |
| 5 | .55E+02 | .55E+00 | .33E+00 |
| 7 | .26E+02 | .34E+00 | .23E+00 |
| 10 | .13E+02 | .23E+00 | .18E+00 |
| 15 | .66E+01 | .16E+00 | .14E+00 |
| 22 | .37E+01 | .12E+00 | .11E+00 |
| 30 | .25E+01 | .98E-01 | .89E-01 |

Table 6. Performances of the Bayes for $\beta = 2$, $n = 5$, $r = 3$ and different $w$ values.

| $RQ[\tilde{x}_R]$ | | | | | $RQ[\tilde{\beta}]$ | | | |
|---|---|---|---|---|---|---|---|---|
| $w$ values | | | | | $w$ values | | | |
| $1.1/\beta$ | $1.4/\beta$ | $1.8/\beta$ | $1/\beta_1+0.1$ | Test # | $1.1/\beta$ | $1.4/\beta$ | $1.8/\beta$ | $1/\beta_1+0.1$ |
| .38E+00 | .28E+00 | .23E+00 | .23E+00 | I | .34E+00 | .25E+00 | .22E+00 | .21E+00 |
| .36E+00 | .11E+01 | .20E+01 | .90E+00 | II | .27E+00 | .40E+00 | .33E+00 | .66E+00 |
| .40E+00 | .46E+00 | .55E+00 | .55E+00 | III | .38E+00 | .46E+00 | .56E+00 | .57E+00 |
| .76E+00 | .72E+00 | .66E+00 | .75E+00 | IV | .78E+00 | .75E+00 | .72E+00 | .86E+00 |
| .84E+00 | .15E+01 | .23E+01 | .12E+01 | V | .91E+00 | .12E+01 | .11E+01 | .40E+00 |
| .76E+00 | .71E+00 | .64E+00 | .55E+00 | VI | .78E+00 | .73E+00 | .67E+00 | .54E+00 |
| .61E+00 | .47E+00 | .39E+00 | .29E+00 | VII | .75E+00 | .57E+00 | .49E+00 | .44E+00 |
| .22E+00 | .72E+00 | .16E+01 | .15E+01 | VIII | .40E+00 | .41E+00 | .50E+00 | .10E+01 |
| .80E+00 | .79E+00 | .81E+00 | .86E+00 | IX | .10E+01 | .10E+01 | .10E+01 | .11E+01 |

Table 6b. Performances of the MLE for $\beta = 2$ and type II censoring.

| $n$ | $r$ | $RQ[\hat{x}_R]$ | $RQ[\hat{\beta}]$ | $DS[\bar{\beta}]$ |
|---|---|---|---|---|
| 5 | 3 | .17E+01 | .12E+02 | .43E+01 |
| 10 | 4 | .11E+01 | .36E+01 | .17E+01 |
| 10 | 6 | .98E+00 | .16E+01 | .10E+01 |
| 20 | 8 | .68E+00 | .12E+01 | .82E+00 |
| 20 | 12 | .60E+00 | .77E+00 | .60E+00 |
| 40 | 16 | .44E+00 | .63E+00 | .52E+00 |
| 40 | 24 | .38E+00 | .45E+00 | .39E+00 |

Table 7. Performances of the Bayes for $\beta = 1$, $n = 5$, $r = 3$ and different $w$ values.

| $RQ[\tilde{x}_R]$ | | | | | $RQ[\tilde{\beta}]$ | | | |
|---|---|---|---|---|---|---|---|---|
| $w$ values | | | | | $w$ values | | | |
| $1.1/\beta$ | $1.4/\beta$ | $1.8/\beta$ | $1/\beta_1+0.1$ | Test # | $1.1/\beta$ | $1.4/\beta$ | $1.8/\beta$ | $1/\beta_1+0.1$ |
| .47E+00 | .40E+00 | .34E+00 | .44E+00 | I | .87E-01 | .70E-01 | .66E-01 | .48E-01 |
| .68E+00 | .14E+01 | .22E+01 | .15E+01 | II | .55E-01 | .87E-01 | .99E-01 | .60E-01 |
| .52E+00 | .53E+00 | .60E+00 | .50E+00 | III | .12E+00 | .14E+00 | .16E+00 | .11E+00 |
| .10E+01 | .94E+00 | .83E+00 | .11E+01 | IV | .13E+00 | .13E+00 | .13E+00 | .15E+00 |
| .12E+01 | .18E+01 | .25E+01 | .15E+01 | V | .15E+00 | .16E+00 | .17E+00 | .15E+00 |
| .10E+01 | .88E+00 | .74E+00 | .97E+00 | VI | .13E+00 | .12E+00 | .11E+00 | .13E+00 |
| .55E+00 | .47E+00 | .41E+00 | .48E+00 | VII | .17E+00 | .15E+00 | .15E+00 | .14E+00 |
| .34E+00 | .85E+00 | .17E+01 | .57E+00 | VIII | .13E+00 | .12E+00 | .12E+00 | .16E+00 |
| .61E+00 | .64E+00 | .68E+00 | .59E+00 | IX | .19E+00 | .19E+00 | .20E+00 | .16E+00 |

Table 7b. Performances of the MLE for $\beta = 1$ and type II censoring.

| $n$ | $r$ | $RQ[\hat{x}_R]$ | $RQ[\hat{\beta}]$ | $DS[\bar{\beta}]$ |
|---|---|---|---|---|
| 5 | 3 | .93E+01 | .58E+01 | .22E+01 |
| 10 | 4 | .48E+01 | .18E+01 | .87E+00 |
| 10 | 6 | .39E+01 | .82E+00 | .51E+00 |
| 20 | 8 | .22E+01 | .60E+00 | .41E+00 |
| 20 | 12 | .19E+01 | .38E+00 | .30E+00 |
| 40 | 16 | .12E+01 | .32E+00 | .26E+00 |
| 40 | 24 | .99E+00 | .22E+00 | .20E+00 |

Table 8. Performances of the Bayes for $\beta = 0.6$, $n = 5$, $r = 3$ and different $w$ values.

| $RQ[\tilde{x}_R]$ | | | | | $RQ[\tilde{\beta}]$ | | | |
|---|---|---|---|---|---|---|---|---|
| $w$ values | | | | | $w$ values | | | |
| $1.1/\beta$ | $1.4/\beta$ | $1.8/\beta$ | $1/\beta_1+0.1$ | Test # | $1.1/\beta$ | $1.4/\beta$ | $1.8/\beta$ | $1/\beta_1+0.1$ |
| .54E+00 | .46E+00 | .38E+00 | .30E+00 | I | .11E+00 | .81E-01 | .78E-01 | .83E-01 |
| .17E+01 | .31E+01 | .41E+01 | .45E+01 | II | .72E-01 | .10E+00 | .12E+00 | .11E+00 |
| .79E+00 | .79E+00 | .83E+00 | .85E+00 | III | .18E+00 | .16E+00 | .16E+00 | .16E+00 |
| .32E+01 | .26E+01 | .20E+01 | .32E+01 | IV | .11E+00 | .11E+00 | .11E+00 | .13E+00 |
| .40E+01 | .46E+01 | .51E+01 | .49E+01 | V | .13E+00 | .15E+00 | .16E+00 | .15E+00 |
| .30E+01 | .21E+01 | .14E+01 | .26E+01 | VI | .95E-01 | .85E-01 | .74E-01 | .10E+00 |
| .75E+00 | .60E+00 | .50E+00 | .56E+00 | VII | .15E+00 | .13E+00 | .11E+00 | .13E+00 |
| .43E+00 | .94E+00 | .18E+01 | .14E+01 | VIII | .11E+00 | .94E-01 | .89E-01 | .12E+00 |
| .91E+00 | .89E+00 | .88E+00 | .88E+00 | IX | .20E+00 | .18E+00 | .17E+00 | .17E+00 |

Table 8b. Performances of the MLE for $\beta = 0.6$ and type II censoring.

| $n$ | $r$ | $RQ[\hat{x}_R]$ | $RQ[\hat{\beta}]$ | $DS[\bar{\beta}]$ |
|---|---|---|---|---|
| 5 | 3 | .91E+02 | .35E+01 | .13E+01 |
| 10 | 4 | .31E+02 | .11E+01 | .52E+00 |
| 10 | 6 | .23E+02 | .49E+00 | .305E+00 |
| 20 | 8 | .93E+01 | .36E+00 | .25E+00 |
| 20 | 12 | .72E+01 | .23E+00 | .18E+00 |
| 40 | 16 | .35E+01 | .19E+00 | .15E+00 |
| 40 | 24 | .28E+01 | .13E+00 | .12E+00 |

## VI. CONCLUSIONS

On the basis of the result of the Monte Carlo study the following conclusions can be drawn:

- When $\beta = 2$ (Tables 3, 3b, 6 and 6b), to obtain values of $RQ[\hat{x}_R]$ similar to the value of $RQ[\tilde{x}_R]$, obtained with $n = 3$ ($n = 5$, $r = 3$) in the case of good prior information (case I) MLE need samples of size $n = 30$ ($n = 40$, $r = 24$).

- When $\beta = 2$, the performances of Bayes estimators, obtained with $w = 1.1/\beta$, $n = 3$ ($n = 5$, $r = 3$) and very biased set of prior information (e.g., see case V) are better than those obtained by MLE with $n = 7$ ($n = 10$, $r = 6$).

- When $\beta = 1$ (Tables 4, 4b, 7 and 7b), to obtain values of $RQ[\hat{x}_R]$ similar to the value of $RQ[\tilde{x}_R]$, obtained with $n = 3$ ($n = 5$, $r = 3$) in the case of good prior information (case I), MLE need samples of size larger than $n = 30$ ($n = 40$, $r = 24$).

- When $\beta = 1$, performances of Bayes estimators, obtained with $w = 1.1/\beta$, $n = 3$ ($n = 5$, $r = 3$) are never worse than those obtained by MLE for $n = 22$ ($n = 40$, $r = 16$).

- When $\beta = 0.6$ [see Tables 5 and 5b (8 and 8b)], values of $RQ[\hat{x}_R]$ obtained with $n = 30$ ($n = 40$, $r = 16$), are almost an order of magnitude greater than the value of $RQ[\tilde{x}_R]$, obtained with $n = 3$ ($n = 5$, $r = 3$), in the case of good prior information (case I);

- When $\beta = 0.6$, the worst performances (case V) obtained by Bayes estimators, with $w = 1.1/\beta$, $n = 3$ ($n = 5$, $r = 3$), are similar to those obtained by MLE with $n = 22$ ($n = 40$, $r = 16$).

- In the considered experimental conditions an alternative to the proposed Bayes-

ian estimators cannot be represented by classical estimators like the MLE (for a comparison between the MLE and other classical estimators see as an example [16] or [18]). In fact, when the experimental samples are very small, MLE supply estimates worse than the Bayes ones, but very often even worse of the prior empirical estimates anticipated on the basis of elementary technological considerations (e.g., the prior information itself);

- Robustness of the Bayes estimators is elevated and moderately depending upon the value assigned to the hyperparameter $w$. Specifically, to higher (lower) values of $w$ it corresponds a minor (greater) ability of the Bayes estimators to react to biased prior information (e.g., cases II, V, VIII).

- High $w$ values (compared to the size of the available sample) render the prior distributions excessively strong and cause automatic confirmation of the prior information. Therefore, using relatively small values of the hyperparameter $w$ (always smaller than the number of failures, $r$) is the best practice in order to not inhibit the physiological Bayesian learning from the experimental data.

ACKNOWLEDGMENT

Dr. R.A. Evans' editorials on the Bayesian topic, published in these Transactions since 1969, have inspired first author's research. Both authors remember Dr. Evans with profound gratitude.